\documentclass[smallextended,envcountsect]{svjour3} 
\usepackage{latexsym, amsmath, amsfonts, amscd, amssymb, verbatim}
\usepackage{graphicx,graphics}
\usepackage{xcolor}
\smartqed 
\newtheorem{Theorem}{Theorem}[section]
\newtheorem{Remark}{Remark}[section]
\newtheorem{Lemma}{Lemma}[section]
\newtheorem{Example}{Example}[section]
\def\R{\rm I\!R}

\def\gph{\mbox{\rm gph}\,}

\def\Limsup{\mathop{{\rm Lim}\,{\rm sup}}}

\def\tto{\;{\lower 1pt \hbox{$\rightarrow$}}\kern -10pt
	\hbox{\raise 2pt \hbox{$\rightarrow$}}\;}

\def\R{I\!\!R}

\def\ox{\bar{x}}
\def\oy{\bar{y}}
\def\oz{\bar{z}}
\def\ov{\bar{v}}

\def\gph{\mbox{\rm gph}\,}

\begin{document}

\title{Sensitivity Analysis of a Stationary Point Set Map under Total Perturbations. Part 1: Lipschitzian Stability}

%\subtitle{Using  the  LaTex Template}

\author{Duong Thi Kim Huyen  \and  Jen-Chih Yao \and Nguyen Dong Yen}

\institute{Duong Thi Kim Huyen \at Graduate Training Center, Institute of Mathematics, Vietnam Academy of Science and Technology,\\ Hanoi, Vietnam\\
 kimhuyenhy@gmail.com
 \and 
 Jen-Chih Yao \at Center for General Education, China Medical University,\\ Taichung, Taiwan\\
 yaojc@mail.cmu.edu.tw
 \and
	Nguyen Dong Yen, Corresponding author \at Institute of Mathematics, Vietnam Academy of Science and Technology,\\ Hanoi, Vietnam\\ ndyen@math.ac.vn}

\date{{\text{Communicated by ...}}
		\\
		\\
Received: date / Accepted: date}
%The correct dates will be entered by the editor.

\maketitle

%\centerline{(Communicated by ....)}

\smallskip
\begin{abstract}
By applying some theorems of Levy and Mordukhovich (Math Program 99: 311--327, {\color{blue}2004}) and other related results, we estimate the Fr\'echet coderivative and the Mordukhovich coderivative of the stationary point set map of a smooth parametric optimization problem with one smooth functional constraint under total perturbations. From the obtained formulas we derive necessary and sufficient conditions for the local Lipschitz-like property of the stationary point set map. This leads us to new insights into the preceding deep
investigations of Levy and Mordukhovich in the above-cited paper and of Qui (J Optim Theory Appl
161: 398--429, {\color{blue} 2014}; J Glob Optim 65: 615--635, {\color{blue} 2016}).
\end{abstract}



%All acknowledgements should be placed in the back of the paper after Conclusions..

\setcounter{equation}{0}
\section{Introduction}
Appeared at the early stage of optimization theory, smooth programming problems continue to attract common attention of the optimization community due to their importance and beauty. Polynomial optimization problems, including nonconvex quadratic programs, are typical examples of such problems.

\smallskip
The present paper investigates the Lipschitz-like property and the Robinson stability of the stationary point set map of a smooth parametric optimization problem with one smooth functional constraint under total perturbations. The aim is achieved by using some theorems  of Levy and Mordukhovich \cite{LeMo04} and other related results from \cite{Yen_Yao_NA_2009}, \cite{LeeYen11AA}, and \cite{Qui_JOTA2014}.

\smallskip
Introduced by Aubin \cite[p.~98]{Aubin_1984} under the name pseudo-Lipschitz property, the \textit{local Lipschitz-like property} of multifunctions is a fundamental concept in stability and sensitivity analysis of optimization and equilibrium problems. It is equivalent to the \textit{classical metric regularity} of the inverse map (see \cite{Borwein_Zhuang_1988,Penot_1989} and \cite{B-M_TAMS_1993}). The local Lipschitz-like property guarantees the local convergence of some variants of Newton's method for generalized equations \cite{Célia_Alain_2005,Dontchev_1996,DR_2009}. In particular, from \cite[Theorem~6C.1, p.~328]{DR_2009} it follows that, if a mild approximation condition is satisfied and the solution map under right-hand-side perturbations is locally Lipschitz-like around a point in question, then there exists an iterative sequence $Q$-linearly converging to the solution. Moreover, as shown by Dontchev \cite[Theorem~1]{Dontchev_1996}, the Newton method applied to a generalized equation in a Banach space is locally convergent uniformly in the canonical parameter if and only if a certain map is locally Lipschitz-like around the reference point. The author also proved (see \cite[Theorem~2]{Dontchev_1996}) that the latter property implies the uniform $Q$-quadratic convergence, provided that the derivative of the base map is locally Lipschitz.

\smallskip
The \textit{Robinson stability} of an implicit multifunction, which has been called the metric regularity in the sense of Robinson by several authors, was introduced by Robinson \cite{Robinson_1976}. This property is a kind of uniform local error bounds and has numerous applications. Recently, Gfrerer and Mordukhovich \cite{Gfrerer_Mor_SIOPT2016} have given first-order and second-order sufficient conditions for this stability property of a parametric constraint system and put it in the relationships with other properties, such as the classical metric regularity and the local Lipschitz-like property.

\smallskip
 The  coderivative analysis of composite constraint functions of Levy and Mordukhovich \cite{LeMo04} is based on the rich generalized differentiation calculus in \cite[Chapter~10]{Rock_Wets_1998}. Among other things, it uses the properties of  amenable functions and strongly amenable functions, and the extended chain rule for subdifferentials \cite[Theorem~10.49]{Rock_Wets_1998}. The analysis allows us to derive sharp upper estimates for the Mordukhovich coderivative of the stationary point set map, where the limiting second-order subdifferential is used.

\smallskip
 To get lower estimates for the Fr\'echet and the Mordukhovich coderivatives of the stationary point set map, we combine the lower estimates of Lee and Yen \cite{LeeYen11AA} with some results of Qui \cite{Qui_JOTA2014,Qui_JoGO_2016}.
 
\smallskip
 With the above upper and lower coderivative estimates, we can use the \textit{Mordukhovich criterion} for the local Lipschitz-like property of locally closed multifunctions to obtain both \textit{necessary and sufficient conditions} for this property. Here, we do not need an additional technical assumption of Qui. Besides, by invoking a result of \cite{Yen_Yao_NA_2009}, we are able to show that these sufficient conditions also guarantee the Robinson stability of the stationary point set map.
 
\smallskip
Our conditions are easy to verify and can be effectively applied to nonconvex quadratic programming under a possibly nonconvex quadratic constraint. The results on quadratic programming in this paper extend the preceding ones of Lee and Yen \cite{LeeYen14NA} and Qui and Yen \cite{QuiYen_SIOPT2014} to a broader class of quadratic programs.

\smallskip
Optimization problems under total perturbations have been studied in \cite{MNg_2014,MNg_2016,MNR_2015,MRS_2013} by different approaches and concepts. But, our results are very different from those of the cited works.

\smallskip
Solution stability of variational inequalities on fixed or linearly perturbed polyhedral convex sets, which is closely related to that of optimization problems under linear constraints, has been investigated intensively; see \cite{Ban_Song_2016,DR_1996,Henrion_Mordukhovich_Nam_2010,Lu_Robinson_2008,Nam_2010,Qui_AMV_2011,Qui_NA_2011,Qui_JMAA_2011,Robinson_2007,Trang_2016,Yao_Yen_AMV_2009,Yao_Yen_PJM_2009}, and the references therein. The case of nonlinearly perturbed polyhedral convex sets has been considered in \cite{Qui_JOTA_2012}.

\smallskip
It is well known that \textit{calmness} is weaker than the local Lipschitz-like property.  Calmness of the stationary point set map of a general parametric optimization problem has been considered, e.g., in \cite[Sect.~4]{Guo_Lin_Ye_SIOPT2012}.

\smallskip
The paper is divided into two parts. In this part, Sect. {\color{blue}2} recalls some basic concepts from variational analysis, formulates the problem studied herein, and presents a series of auxiliary results in a unified form. Sections {\color{blue}3} and~{\color{blue}4} present new results on smooth parametric optimization problem with one smooth functional constraint under total perturbations. Namely, sensitivity analysis of the stationary point set at the interior points (resp., at the boundary points) of the progamming variable-parameter domain is given in Sect.~{\color{blue}3} (resp., in Sect.~{\color{blue}4}). 

\smallskip
In Part 2,  sufficient conditions for the Robinson stability of the stationary point set map will be established. This allows us to revisit and extend several stability theorems in indefinite quadratic programming. A comparison of our results with the ones which can be obtained via another approach will be also given. 

\section{Preliminaries}
The scalar product and the norm in a finite-dimensional Euclidean space are denoted respectively by $\langle
\cdot,\cdot \rangle$ and $\|\cdot \|$. The symbols $B(x,\rho)$ and $\bar{B}(x,\rho)$ stand for the open (resp., closed) ball centered at $x \in X$ with radius $\rho > 0$. The distance $\displaystyle\inf_{u\in A}\|x-u\|$ from $x \in X$ to a subset $A \subset X$ is denoted by $d(x, A)$.

We now recall several basic concepts from variational analysis \cite{Rock_Wets_1998,B-M06} which will be used intensively later on.

The {\it Fr\'echet normal cone} (also called the {\it prenormal
	cone}, or the {\it regular normal cone}) to a set
$\Omega\subset\R^s$ at $\ov \in \Omega$ is given by
\begin{eqnarray*}\label{Frechet normals} \widehat
	N_{\Omega}(\ov)=\left\{v'\in\R^s\, \mid\,
	\displaystyle\limsup_{v\xrightarrow{\Omega}\ov}\,\displaystyle\frac{\langle
		v',v-\ov \rangle}{\|v-\ov \|}\leq 0\right\},\end{eqnarray*} where
$v\xrightarrow{\Omega}\ov$ means $v\to \ov$ with $v\in\Omega$. By
convention, $\widehat N_{\Omega}(\ov):=\emptyset$ when $\ov \notin \Omega$.
Provided that $\Omega$ is locally closed around $\bar v\in \Omega$, one calls
\begin{eqnarray*}\label{basic normals}\begin{array}{rl}
		N_{\Omega}(\ov)&=\displaystyle\Limsup_{v\to\ov}\widehat
		N_{\Omega}(v)\\
		&:=\big\{v'\in\R^s\, \mid\,  \exists \mbox{ sequences } v_k\to \ov,\ v_k'\rightarrow v',\\
		& \qquad \qquad \qquad \quad \ \mbox {with } v_k'\in \widehat
		N_{\Omega}(v_k)\, \mbox{ for all }\, k=1,2,\dots
		\big\}\end{array}\end{eqnarray*} the {\it Mordukhovich} (or {\it
	limiting/basic}) {\it normal cone} to $\Omega$ at $\ov$.
If $\ov \notin \Omega$, then one puts $N_{\Omega}(\ov)=\emptyset$.

A multifunction $\Phi:\R^n\rightrightarrows  \R^m$ is
said to be {\it locally closed} around a point $\oz=(\ox,\oy)$ from  $\gph\Phi:=\{(x,y)\in \R^n\times\R^m\,\mid\, y\in\Phi(x)\}$ if $\gph\Phi$ is locally closed around $\bar z$. Here, the product space
$\R^{n+m}=\R^n\times\R^m$ is equipped with the topology generated by the sum norm $\|(x,y)\| =\|x\|+\|y\|$.

For any $\bar{z}=(\bar{x},\bar{y})\in \mbox{gph}\,\Phi$,
\begin{eqnarray*}\label{Frechet coderivative}
	\widehat D^*\Phi(\bar{z})(y'):=\big\{x'\in
	\R^n\, \mid\, (x',-y')\in \widehat N_{\mbox{gph}\,\Phi}(\bar{z})\big\}\quad (y'\in \R^m)
\end{eqnarray*} are called the {\it
Fr\' echet coderivative} values of $\Phi$ at $\bar z$. Similarly, the {\it
Mordukhovich coderivative} (limiting coderivative) values of $\Phi$ at $\bar z$ are defined by
\begin{eqnarray*}\label{normal coderivative}D^*\Phi(\bar{z})(y'):=\big\{x'\in
	\R^n\, \mid\, (x',-y')\in N_{\mbox{gph}\,\Phi}(\bar{z})\big\}\quad (y'\in \R^m).\end{eqnarray*}
Thus, $\widehat D^*\Phi(\bar{z})$ and $D^*\Phi(\bar{z})$ are multifuntions from $\R^m$ to $\R^n$. By \cite[Theorem~1.38]{B-M06}, if $\Phi$ is strictly Fr\' echet differentiable at $\bar x$, then
$$\widehat D^*\Phi(\bar{x})(y')=D^*\Phi(\bar{x})(y')
=\{\nabla\Phi(\bar x)^*(y')\}$$
for any $y'\in \R^m$.

Suppose that $X$, $Y$, and $Z$ are finite-dimensional Euclidean spaces. Consider a function $\psi: X\rightarrow \bar{\R}$ with $|\psi(\bar x)|<\infty$. The set
$$\partial\psi(\bar x):=\{x'\in X^*\mid (x',-1)\in N_{{\rm epi}\,\psi}(\bar x,\psi(\bar x))\}$$
is the \textit{Mordukhovich subdifferential} of $\psi$ at $\bar x$. We put $\partial\psi(\bar x)=\emptyset$ if $|\psi(\bar x)|=\infty$. The set
$$\partial^{\infty}\psi(\bar x):=\{x^*\in X^*\mid (x^*,0)\in N_{{\rm epi}\,\psi}(\bar x,\psi(\bar x))\}$$
is the \textit{singular subdifferential} of $\psi$ at $\bar x$. For a set $\Omega\subset X$ and a point $\bar x\in \Omega$, we have
$$N_{\Omega}(\bar x)=\partial\delta_{\Omega}(\bar x)=\partial^{\infty}\delta_{\Omega}(\bar x),$$
where $\delta_{\Omega}(\bar x)$ is the indicator function of $\Omega$; see \cite[Proposition~1.79]{B-M06}. If $\psi$ depends on two variables $x$ and $y$, and $|\psi(\bar x, \bar y)|<\infty$, then $\partial_x\psi(\bar x,\bar y)$ denotes the Mordukhovich subdifferential of $\psi(.,\bar y)$ at $\bar x$. For any $\bar v\in \partial\psi(\bar x)$,
$$\partial^2\psi(\bar x|\bar v)(u):=D^*(\partial\psi)(\bar x|\bar v)(u)\quad (u\in X^{**}=X)$$ is the \textit{limiting second-order subdifferential} (or the generalized Hessian).

A multifunction $G: Y \rightrightarrows X$ is said to be \textit{locally Lipschitz-like} around $(\bar{y}, \bar{x})\in {\rm gph}\,G$ if there exists a constant $\ell > 0$ and neighborhoods $U$ of $\bar{x}$, $V$ of $\bar{y}$ such that
\begin{equation*}\label{Lipschitz-like property} G(y') \cap U \subset G(y) + \ell \|y' - y\| \bar{B}_{X} \quad \forall y, y' \in V, \end{equation*}
where $\bar{B}_{X}$ denotes the closed unit ball in $X$. When $G$ is locally closed around $(\bar{y}, \bar{x})$, the \textit{Mordukhovich criterion} (see \cite{B-M_TAMS_1993}, \cite[Theorem~9.40]{Rock_Wets_1998}, and \cite[Theorem~4.10]{B-M06}) says that $G$ is locally Lipschitz-like around $(\bar{y}, \bar{x})$ if and only if
\begin{equation*}\label{M_criterion}D^{*}G(\bar{y}, \bar{x})(0) = \{0\}. \end{equation*}

For a multifunction $F: X \times Y \rightrightarrows Z$ and a pair $(\bar{x}, \bar{y}) \in X\times Y$ satisfying $0\in F(\bar{x}, \bar{y})$, we say that the \textit{implicit multifunction} $G: Y \rightrightarrows X$ given by
$G(y) = \{x \in X\mid 0 \in F(x, y)\}$ has the \textit{Robinson stability} at $\omega_{0}:= (\bar{x}, \bar{y}, 0)$ if there exist constants $r > 0$, $\gamma > 0$, and neighborhoods $U$ of $\bar{x}$, $V$ of $\bar{y}$ such that
\begin{equation*}\label{Metric_Regularity}d(x, G(y)) \leq rd(0, F(x, y)) \end{equation*}
for any $(x, y)\in U\times V$ with $d(0, F(x, y)) < \gamma$.
Note that the condition $d(0, F(x, y)) < \gamma$ can be omitted if $F$ is \textit{inner semicontinuous} at $(\bar x, \bar y, 0)$; see \cite{Huyen_Yen_2016}. Note that, in some cases, the Robinson stability of $G$ at $(\bar x,\bar y,0)$ implies its local Lipschitz-likeness around $(\bar y,\bar x)$; see, e.g., [8]. For the generalized linear constraint system studied in [11], these properties are equivalent. In the sequel, we will see that the regularity conditions in use guarantee for our stationary point set map to have both properties.

Now, let $f_0$ and $F$ be twice continuously differentiable real-valued functions ($C^2$-functions for brevity) defined on the product $\R^n\times \R^d$ of two Euclidean spaces. For every $w\in\R^d$, we consider the parametric optimization problem
\begin{equation*}\label{Optim_Prob}
(P_w)\quad\ \, {\rm Minimize}\ \, f_0(x,w) \ \; {\rm subject\ to}\ x\in \R^n\ \, {\rm and}\ \, F(x,w)\leq 0.
\end{equation*} The constraint set of $(P_w)$ is $C(w):=\{x\in\R^n\mid F(x,w)\leq 0\}$. The stationary point set of $(P_w)$ is defined by
\begin{equation}\label{KKT_point-set0}
S(w)=\{x\in \R^n\mid 0\in  \nabla_{x}f_0(x,w)+N_{C(w)}(x)\}.
\end{equation}
When $w$ varies on $\R^d$, one has a multifunction  $S:\R^d\rightrightarrows\R^n$ with $S(w)$ being calculated by~\eqref{KKT_point-set0}. Setting $f(x,w)=g(F(x,w))=(g\circ F)(x,w)$, where $g(y)=\delta_{\R_-}(y)$, i.e., $g(y)=0$ for $y\in (-\infty,0]$ and $g(y)=+\infty$ for $y>0$, we can rewrite \eqref{KKT_point-set0} as
\begin{equation}\label{KKT_point-set}
S(w)=\{x\in \R^n\mid 0\in  \nabla_{x}f_0(x,w)+\partial_{x}f(x,w)\}.
\end{equation}

Fix a vector $w=\bar w\in \R^d$ and suppose that $\bar x\in S(\bar w)$. Since $(P_{\bar w})$ has a single smooth inequality constraint, the Mangasarian-Fromovitz Constraint Qualification is fulfilled at $\bar x\in C(\bar w)$ if and only if
\begin{equation}\label{MFCQ_condition}
\text{If}\ F(\bar x,\bar w)=0,\ \text{then}\ \nabla_xF(\bar x,\bar w)\neq 0.\quad\quad\tag{\textbf{MFCQ}}
\end{equation}
In what follows, we assume that \textbf{(MFCQ)} is valid.
To study the stability of the stationary point set map $S$ around the $(\bar w,\bar x)$ in ${\rm gph}\,S$, we compute the Mordukhovich and the Fr\' echet coderivatives of the partial subdifferential map $\partial_{x}f:\R^n\times \R^d\rightrightarrows\R^n$. In general, there is no explicit formula for the coderivatives of such maps. However, the results of \cite{LeMo04} provide us with some tools which allow us to estimate the coderivative value $D^*S(\bar w|\bar x)(x')$ for every $x'\in\R^n$.

The fulfillment of MFCQ at $(\bar x,\bar w)$ implies that $g(x,w)=g(F(x,w))$ is a strongly amenable in $x$ at $\bar x$ with compatible parameterization in $w$ at $\bar w$. Then, by \cite[Theorem~10.49]{Rock_Wets_1998}, for $(x,w)$ near $(\bar x,\bar w)$, we have
\begin{equation}\label{Implication_of_s_amenability1}
\partial f(x,w)=\nabla F(x,w)^*(\partial g(F(x,w)))
\end{equation}
and
\begin{equation}\label{Implication_of_s_amenability2}
\partial_xf(x,w)=\nabla_xF(x,w)^*(\partial g(F(x,w)));
\end{equation}
see \cite[formulas~(14) and (15)]{LeMo04}.

In order to estimate the limiting second-order subdifferential of $f$, we need the following result.

\begin{Lemma}\label{Theorem3.1_LeMo04} {\rm (see \cite[Theorem 3.1]{LeMo04})} 
	Suppose that $\bar v\in \partial f(\bar x,\bar w)$. Then, for
	any $v'\in \R^n\times \R^d$,
	\begin{equation*}
	\begin{array}{rl}
	&\partial^{2}f((\bar x,\bar w)|\bar v)(v')\\
	& \subset\ \displaystyle\bigcup_
	{\begin{subarray}{c} \bar y\in \partial g(F(\bar x,\bar w))\; \text{with}\\ \nabla F(\bar x,\bar w)^*\bar y=\bar v \end{subarray}}
	\bigg(\nabla^{2}(\bar {y}\cdot F)(\bar x,\bar w)v'
	+D^*(\partial g\circ F)(\bar x,\bar w)|\bar y)(\nabla F(\bar x,\bar w)v')\bigg),
	\end{array}
	\end{equation*}
	where the function $\bar {y}\cdot F:\R^{n+d}\to \R$ is defined by $(\bar {y}\cdot F)(x,w):=\bar y F(x,w)$. If, in addition, at every $\bar y\in \partial g(F(\bar x,\bar w))$ with $\nabla F(\bar x,\bar w)^*\bar y=\bar v$, one has the second-order constraint qualification
	\begin{equation}\label{Condi17_LeMor}
	\partial^{2}g(F(\bar x,\bar w)|\bar y)(0)\cap {\rm ker}\nabla F(\bar x,\bar w)^*=\{0\},
	\end{equation} then the estimate above for the second-order subdifferential can be refined by replacing the coderivative of the multifunction $\partial g\circ F$ via the inclusion
	\begin{equation*}\label{Upper_esti_coder_compositefunc}
	D^*(\partial g\circ F)((\bar x,\bar w)|\bar y)(\nabla F(\bar x,\bar w)v')\subset \nabla F(\bar x,\bar w)^* \partial^{2}g(F(\bar x,\bar w)|\bar y)(\nabla F(\bar x,\bar w)v').
	\end{equation*}
\end{Lemma}

In our problem $(P_w)$, condition \eqref{Condi17_LeMor} can be omitted. Indeed, $\bar y\in\partial g(F(\bar x,\bar w))$ if and only if $\bar y \in N_{\R_{-}}(F(\bar x,\bar w))$. Hence, $\bar y\geq 0$. Clearly, $${\rm gph}\,\partial g=(\R_{-}\times \{0\})\cup (\{0\}\times \R_{+}).$$
If $F(\bar x,\bar w)<0$, then $\bar y=0$ and $N_{{\rm gph}\,\partial g}(F(\bar x,\bar w),\bar y)=\{0\}\times \R.$
It follows that
\begin{eqnarray*}\begin{array}{rl} \partial^2g(F(\bar x,\bar w)|\,\bar y)(0)& =D^*(\partial g(F(\bar x,\bar w)|\bar y))(0)\\
		& =\{u'\in \R\mid (u',0)\in N_{{\rm gph}\,\partial g}(F(\bar x,\bar w),\bar y)\}=\{0\}.
	\end{array}
\end{eqnarray*}
So \eqref{Condi17_LeMor} is satisfied. If $F(\bar x,\bar w)=0$, then \textbf{(MFCQ)} implies $\nabla F(\bar x,\bar w)\neq 0$. Hence the linear operator $\nabla F(\bar x,\bar w):\R^n\times\R^n\to\R$ is surjective. Thus ${\rm ker}\,\nabla F(\bar x,\bar w)^*=\{0\}$ by~\cite[Lemma~1.18]{B-M06}, and we see that \eqref{Condi17_LeMor} is fulfilled. Therefore, applied to $(P_w)$, Lemma~\ref{Theorem3.1_LeMo04} can be reformulated as follows: \textit{For any} $\bar v\in\partial f(\bar x,\bar w)$ \textit{and} $v'\in \R^n\times \R^d$,
\begin{equation}\label{Upper_esti_2nd_order_subdifferential_of_f}
\partial^{2}f((\bar x,\bar w)|\bar v)(v')\subset\displaystyle\bigcup_
{\begin{subarray}{c} \bar y\in \partial g(F(\bar x,\bar w))\;\text{with}\\ \nabla F(\bar x,\bar w)^*\bar y=\bar v \end{subarray}}
\Big(\nabla^{2}(\bar {y}\cdot F)(\bar x,\bar w)v'
+\Omega_1(\bar y,v')\Big),
\end{equation} \textit{where}
\begin{equation*}\label{Omega_1}
\Omega_1(\bar y,v'):=\nabla F(\bar x,\bar w)^* \partial^{2}g(F(\bar x,\bar w)|\bar y)(\nabla F(\bar x,\bar w)v').
\end{equation*}

\begin{Remark} {\rm Concerning the paper \cite{MR_2012}, 	
	observe that the set $\partial^{2}f((\bar x,\bar w)|\bar v)(v')$ in  formula \eqref{Upper_esti_2nd_order_subdifferential_of_f} is analogous to the set $\widetilde\varphi^2_x(\bar x,\bar w,\bar y)(u)$ (a value of the extended partial second-order subdifferential) in formula~(3.4) of that work. A careful checking shows that equality (3.4) of \cite{MR_2012} implies the upper estimate \eqref{Upper_esti_2nd_order_subdifferential_of_f}.}
\end{Remark}

In what follows, for any $\bar v=(\bar v_x,\bar v_w)\in\R^n\times \R^d$, we put ${\rm proj}_1\bar v=\bar v_x$. The upper estimation for the coderivative values of the stationary point set map $S$ given by Levy and Mordukhovich \cite{LeMo04} requires the following \textit{regularity condition}: \textit{For any} $v'_1\in\R^n$,
\begin{equation}\label{Condi11_LeMor}
0\in \nabla^{2}f_0(\bar x,\bar w)^*(v'_1,0)+\displaystyle\bigcup_{\begin{subarray}{c} \bar v\in \partial f(\bar x,\bar w)\;\text{with}\\ {\rm proj}_1\bar v=-\nabla_{x}f_0(\bar x,\bar w)\end{subarray}}\partial^2f((\bar x,\bar w)|\bar v)(v'_1,0)\ \; \Longrightarrow\ \; v'_1=0
\end{equation}
(see \cite[formula~(11)]{LeMo04}).
For our problem $(P_w)$, by the assumption \textbf{(MFCQ)} and formula~\eqref{Implication_of_s_amenability1}, we have $\partial f(\bar x,\bar w)=\nabla F(\bar x,\bar w)^*(\partial g(\bar x,\bar w))$. In addition, it is easy to show that, for every $\bar y\in\partial g(\bar x,\bar w)$, ${\rm proj}_1\left(\nabla F(\bar x,\bar w)^*\bar y\right)=\nabla_xF(\bar x,\bar w)^*\bar y$. Hence
\begin{equation*}\begin{array}{rl}
& \displaystyle\bigcup_{\begin{subarray}{c} \bar v\in \partial f(\bar x,\bar w)\;\text{with}\\ {\rm proj}_1\bar v=-\nabla_{x}f_0(\bar x,\bar w)\end{subarray}}\partial^2f((\bar x,\bar w)|\bar v)(v'_1,0)\\
&=\displaystyle\bigcup_{\begin{subarray}{c} \bar y\in \partial g(F(\bar x,\bar w))\;\text{with}\\\nabla_{x}F(\bar x,\bar w)^*\bar y=-\nabla_{x}f_0(\bar x,\bar w)\end{subarray}}\partial^2f((\bar x,\bar w)|\nabla F(\bar x,\bar w)^*\bar y)(v'_1,0).
\end{array}\end{equation*}
So \eqref{Condi11_LeMor} is equivalent to the following condition:
\begin{equation}\label{C0_condition}
0\in \nabla^{2}f_0(\bar x,\bar w)^*(v'_1,0)+\Omega_2(v'_1)\ \; \Longrightarrow\ \; v'_1=0,\quad\quad\tag{\textbf{C0}}
\end{equation}
where \begin{equation}\label{Omega_2} \Omega_2(v'_1):=\displaystyle\bigcup_{\begin{subarray}{c} \bar y\in \partial g(F(\bar x,\bar w))\;\text{with}\\\nabla_{x}F(\bar x,\bar w)^*\bar y=-\nabla_{x}f_0(\bar x,\bar w)\end{subarray}}\partial^2f((\bar x,\bar w)|\nabla F(\bar x,\bar w)^*\bar y)(v'_1,0).\end{equation}

The next result from \cite{LeMo04} provides us with an upper estimation for the values of the coderivative map  $D^*S(\bar w|\bar x): \R^n\rightrightarrows\R^d$.

\begin{Lemma}[{see \cite[Corollary~3.1]{LeMo04}}]\label{Corollary3.1_LeMo04}
	If the regularity condition \eqref{C0_condition} holds then, for each $x'\in\R^n$, the coderivative value $D^*S(\bar w|\bar x)(x')$ is contained in the set of $w'\in \R^d$ for which there exists a vector $v'_1\in \R^n$ with
	\begin{equation*}
	(-x',w')-\nabla^{2}f_0(\bar x,\bar w)^*(v'_1,0)\in \Omega_2(v'_1).
	\end{equation*}
\end{Lemma}

Although it is rather difficult to compute the set $\Omega_2(v'_1)$, we can still estimate it by using \eqref{Upper_esti_2nd_order_subdifferential_of_f}.

\textit{Upper estimates} for the limiting coderivative values of $S$ can be derived from a result of Levy and Mordukhovich \cite[Theorem~2.1]{LeMo04}. But, a constraint qualification must be imposed to have these estimates (see \cite[p.~1020]{LeeYen11AA} for details). Interestingly, due to a result of Lee and Yen \cite[Theorem~3.4]{LeeYen11AA}, sharp \textit{lower estimates} for the Fr\' echet coderivative values of $S$ can be given without any condition. Put $G(x,w)=\nabla_xf_0(x,w)$ and $M(x,w)=\partial_{x}f(x,w)$. Then,
\begin{equation}\label{Sum_for_S} S(w)=\{x\in\R^n\mid 0\in G(x,w)+M(x,w)\}.\end{equation}
Since $\bar x\in S(\bar w)$, $\bar\tau:=(\bar x,\bar w,-\nabla_xf_0(\bar x,\bar w))$ belongs to ${\rm gph}\,M$. Note that ${\rm gph}\,M$ is locally closed around $\bar\tau$. The following result combines the lower estimates with the upper estimates mentioned above.

\begin{Lemma}\label{Combined_estimates} {\rm (see \cite[Theorem~3.4]{LeeYen11AA})}
	The lower estimates
	$$\widehat\Gamma(x')\subset\widehat D^*S(\bar w|\bar x)(x')\subset D^*S(\bar w\mid\bar x)(x'),$$
	where
	\begin{equation}\label{Widehat_Gamma(x')}
	\widehat\Gamma(x'):=\displaystyle\bigcup_{v'_1\in \R^n}\left\{w'\in \R^d\mid (-x',w')\in \nabla G(\bar x,\bar w)^*v'_1+\widehat D^*M(\bar\tau)(v'_1)\right\},
	\end{equation} hold for any $x'\in\R^n$. If the constraint qualification
	\begin{equation}\label{C1_condition}
	0\in \nabla G(\bar x,\bar w)^*v'_1+D^*M(\bar\tau)(v'_1)\ \; \Longrightarrow\ \;v'_1=0
	\quad\quad\tag{\textbf{C1}}
	\end{equation}
	is satisfied, then the upper estimate
	$$D^*S(\bar w|\bar x)(x')\subset\Gamma(x'),$$ where
	\begin{equation*}\label{Gamma(x')}
	\Gamma(x'):=\displaystyle\bigcup_{v'_1\in \R^n}\left\{w'\in \R^d\mid (-x',w')\in \nabla G(\bar x,\bar w)^*v'_1+ D^*M(\bar\tau)(v'_1)\right\},
	\end{equation*} is valid for any $x'\in\R^n$. If, in addition, M is graphically regular at $\bar\tau$, then
	$$\widehat\Gamma(x')=\widehat D^*S(\bar w|\bar x)(x')= D^*S(\bar w|\bar x)(x')=\Gamma(x').$$
\end{Lemma}

From Lemma \ref{Combined_estimates}, for any $x'\in\R^n$, $\widehat\Gamma(x')\subset\widehat D^*S(\bar w|\bar x)(x')$.
This implies that $\widehat\Gamma(0)\subset\widehat D^*S(\bar w|\bar x)(0)\subset D^*S(\bar w|\bar x)(0)$. If we put $\widetilde{M}(x,w)=G(x,w)+M(x,w)$, then by the Fr\'echet coderivative sum rule with equalities \cite[Theorem~1.62]{B-M06},
\begin{equation*}\label{Sum_rule_Fr1}
\widehat D^*\widetilde{M}(\omega_0)(v'_1)=\nabla G(\bar x,\bar w)^*v'_1+\widehat D^*M(\bar\tau)(v'_1)
\end{equation*}
for any $v'_1\in\R^n$, where $\omega_0:=(\bar x,\bar w,0)\in {\rm gph}\,\widetilde{M}$.  Therefore, we can write 	$$\widehat\Gamma(x')=\displaystyle\bigcup_{v'_1\in \R^n}\left\{w'\in \R^d\mid (-x',w')\in\widehat D^*\widetilde{M}(\omega_0)(v'_1)\right\}.$$
Note that $0\in \widehat\Gamma(0)$. According to the Mordukhovich criterion, if $S$ is locally Lipschitz-like around $(\bar w,\bar x)$, then $D^*S(\bar w|\bar x)(0)=\{0\}$ and $\widehat\Gamma(0)=\{0\}$ as a result. In addition, if the constraint qualification \eqref{C1_condition} is fulfilled, then Lemma~\ref{Combined_estimates} yields $D^*S(\bar w|\bar x)(x')\subset\Gamma(x')$ for any $x'\in\R^n$.
In particular, $D^*S(\bar w|\bar x)(0)\subset\Gamma(0)$. Hence, if \eqref{C1_condition} is valid and $\Gamma(0)=\{0\}$, then $$D^*S(\bar w|\bar x)(0)=\{0\}.$$ So, due to the Mordukhovich criterion, $S$ is locally Lipschitz-like around $(\bar w,\bar x)$. This idea has been presented in \cite{LeeYen11AA} and we will follow it throughout this paper.

Put $\mathcal{D}=\{(x,w)\in \R^n\times\R^d\mid F(x,w)\leq 0\}$. If $F(\bar x,\bar w)<0$, then $(\bar x, \bar w)$ is an \textit{interior point} of~$\mathcal{D}$. If  $F(\bar x,\bar w)=0$, then $(\bar x, \bar w)$ is a \textit{boundary point} of $\mathcal{D}$. In the next two sections, we will consider separately these two possibilities of the reference point $(\bar x, \bar w)$. Remind that $\bar w\in \R^d$ and $\bar x\in S(\bar w)$ are fixed and all the notations of this section are kept unchanged.

\section{Interior points}

Suppose that $F(\bar x,\bar w)<0$, i.e.,  $(\bar x, \bar w)$ is an \textit{interior point} of~$\mathcal{D}$. A point $\bar x$ belongs to $S(\bar w)$ if and only if
$$0\in  \nabla_{x}f_0(\bar x,\bar w)+N_{C(\bar w)}(\bar x),$$ where $C(\bar w)=\{x\in\R^n\mid F(x,\bar w)\leq 0\}$. Since $F(\bar x,\bar w)<0$, the continuity of $F(.,\bar w)$ implies that $\bar x\in{\rm int}\,C(\bar w)$. This yields $N_{C(\bar w)}(\bar x)=\{0\}$. Thus, $\bar x\in S(\bar w)$ if and only if $\nabla_{x}f_0(\bar x,\bar w)=0$.

The inequality $F(\bar x,\bar w)<0$ implies that $\partial g(F(\bar x,\bar w))=\{0\}$. So, $\bar y=0$ is the unique element of $\partial g(F(\bar x,\bar w))$. Since ${\rm gph}\,\partial g=(\R_{-}\times \{0\})\cup (\{0\}\times\R_{+})$,
$$N_{{\rm gph}\,\partial g}(F(\bar x,\bar w),\bar y)=\{0\}\times \R.$$
Hence, for any $v'\in\R^n\times\R^d$,
\begin{eqnarray*}\begin{array}{rl} & \partial^2g(F(\bar x,\bar w)|\,\bar y)(\nabla F(\bar x,\bar w)v')\\ 
		& =D^*(\partial g(F(\bar x,\bar w)|\bar y))(\nabla F(\bar x,\bar w)v')\\
		& =\{u'\in \R\mid (u',-\nabla F(\bar x,\bar w)v')\in N_{{\rm gph}\,\partial g}(F(\bar x,\bar w),\bar y)\}=\{0\}.		
	\end{array}\end{eqnarray*}
	Therefore, $\Omega_1(\bar y,v')=\{0\}$ for any $v'\in\R^n\times\R^d$. So, from \eqref{Upper_esti_2nd_order_subdifferential_of_f} it follows that
	$$\partial^{2}f((\bar x,\bar w)|\bar v)(v')\subset \{0\}$$
	for any $\bar v\in\partial f(\bar x,\bar w)$ and  $v'\in \R^n\times \R^d$. Since $\nabla_{x}f(\bar x,\bar w)=0$, invoking \eqref{Implication_of_s_amenability1} and the fact that $\partial g(F(\bar x,\bar w))=\{0\}$, we get
	$$\Omega_2(v'_1)=\partial^2f((\bar x,\bar w)|0)(v'_1,0)\subset\{0\}$$ for any $v'_1\in\R^n$.
	Then, condition \eqref{C0_condition} is fulfilled if
	$$0\in \nabla^{2}f_0(\bar x,\bar w)^*(v'_1,0)\ \;\Longrightarrow\ \;v'_1=0.$$
	By the symmetry of $\nabla^{2}_{xx}f_0(\bar x,\bar w)$ and the equality $\nabla^{2}_{xw}f_0(\bar x,\bar w)^T=\nabla^{2}_{wx}f_0(\bar x,\bar w)$, this is equivalent to
	$$\begin{cases}
	\nabla^{2}_{xx}f_0(\bar x,\bar w)v'_1=0\\
	\nabla^{2}_{wx}f_0(\bar x,\bar w)v'_1=0
	\end{cases}\ \;\Longrightarrow\ \;v'_1=0.$$
	Clearly, the latter means that
	\begin{equation}\label{Condi_tildeC1_S}
	{\rm ker}\,\nabla^{2}_{xx}f_0(\bar x,\bar w)\cap {\rm ker}\,\nabla^{2}_{wx}f_0(\bar x,\bar w)=\{0\}.
	\end{equation}
	We now suppose that condition \eqref{Condi_tildeC1_S}, which guarantees the validity of \eqref{C0_condition}, is satisfied. Then, by Lemma~\ref{Corollary3.1_LeMo04},
	$$D^*S(\bar w|\bar x)(x')\subset \displaystyle\bigcup_{v'_1\in \R^n}\{w'\in \R^d\mid (-x',w')-\nabla^{2}f_0(\bar x,\bar w)^*(v'_1,0)\in \Omega_2(v'_1)\},$$
	for any $x'\in \R^n$. Since $\Omega_2(v'_1)\subset\{0\}$, we have
	\begin{equation*}\label{Gamma_1(x')}
	\begin{array}{rl}
	& D^*S(\bar w|\bar x)(x')\\
	& \subset\Gamma_1(x') := \displaystyle\bigcup_{v'_1\in \R^n}\left\{w'\in \R^d\mid \nabla^{2}_{xx}f_0(\bar x,\bar w)v'_1=-x',\; w'=\nabla^{2}_{wx}f_0(\bar x,\bar w)v'_1\right\}.
	\end{array}
	\end{equation*}
	Note that $0\in D^*S(\bar w|\bar x)(0)$. So, if $\Gamma_1(0)=\{0\}$, then $D^*S(\bar w|\bar x)(0)=\{0\}$; as a result, $S$ is locally Lipschitz-like around $(\bar w,\bar x)$ by the Mordukhovich criterion. We have $\Gamma_1(0)=\{0\}$ if and only if
	$$\left\{w'\in \R^d\mid \exists v'_1\in \R^n\ {\rm with}\ \nabla^{2}_{xx}f_0(\bar x,\bar w)v'_1=0,\; w'=\nabla^{2}_{wx}f_0(\bar x,\bar w)v'_1\right\}=\{0\}.$$
	This can be rewritten equivalently as
	\begin{equation}\label{Condi_tildeC2}
	{\rm ker}\,\nabla^{2}_{xx}f_0(\bar x,\bar w)\subset{\rm ker}\,\nabla^{2}_{wx}f_0(\bar x,\bar w).
	\end{equation}
	It is easy to show that \eqref{Condi_tildeC2} and \eqref{Condi_tildeC1_S} hold simultaneously if and only if
	\begin{equation}\label{Aubin_Suffi_Condi1}
	{\rm ker}\,\nabla^{2}_{xx}f_0(\bar x,\bar w)=\{0\}.
	\end{equation} In particular, \eqref{Aubin_Suffi_Condi1} is a \textit{sufficient condition} for $S$ being locally Lipschitz-like around $(\bar w,\bar x)$.
	
	\begin{Example} {\rm 
		Consider the problem $(P_w)$ with $f_0(x,w)=\frac{1}{2}x^TDx+c^Tx$ and $F(x,w)=\|x\|^2-\rho^2$, where $w=(D,c,\rho)$ with $D$ being a  $n\times n$ symmetric matrix, $c\in\R^n$, and $\rho>0$. Suppose that $\bar x\in S(\bar w)$ with $\bar w:=(\bar D,\bar c,\bar\rho)$ and $\|\bar x\|<\bar \rho$. If $\det\bar D\neq 0$, then $S$ is locally Lipschitz-like around $(\bar w,\bar x)$ because \eqref{Aubin_Suffi_Condi1} is satisfied.}
	\end{Example}
	
	Having the sufficient condition \eqref{Aubin_Suffi_Condi1} for the  Lipschitz-likeness of $S$ around $(\bar w,\bar x)$, we want to find a \textit{necessary condition} for this property. We know that if $S$ is locally Lipschitz-like around $(\bar w,\bar x)$, then $\widehat\Gamma(0)=\{0\}$, where the sets have been defined in \eqref{Widehat_Gamma(x')}. Since $F(\bar x,\bar w)<0$ and $F$ is continuous, there exit neighborhoods $U$ of $\bar x$ and $W$ of $\bar w$ such that $F(x,w)<0$ for any $(x,w)\in U\times W$. It follows that, for each $w\in W$, the inclusion $x\in {\rm int}\,C(w)$ holds for every $x\in U$. Hence, for each $w\in W$, $N_{C(w)}(x)=\{0\}$ for all $x\in U$. This means that $M(x,w)=\{0\}$ for $(x,w)$ in a neighborhood of $(\bar x,\bar w)$. Therefore, from \eqref{KKT_point-set0}, 
\begin{equation}\label{new_S} S(w)\cap U=\{x\in\R^n\mid \nabla_xf_0(x,w)=0\}\cap U\quad (\forall w\in W).\end{equation} Since $\widehat D^*M(\bar\tau)(v'_1)=\{0\}$ for any $v'_1\in \R^n$, one has
	$$
	\begin{array}{rl}
	\widehat\Gamma(x') & =
	\displaystyle\bigcup_{v'_1\in \R^n}\left\{w'\in \R^d\mid (-x',w')\in \nabla G(\bar x,\bar w)^*v'_1\right\}\\
	& =\displaystyle\bigcup_{v'_1\in \R^n}\left\{w'\in \R^d\mid (-x',w')\in (\nabla^2_{xx} f_0(\bar x,\bar w)v'_1,\nabla^2_{wx} f_0(\bar x,\bar w)v'_1)\right\}.
	\end{array}$$
	It follows that
	$$\widehat\Gamma(0)=\displaystyle\bigcup_{v'_1\in \R^n}\left\{w'\in \R^d\mid (0,w')\in (\nabla^2_{xx} f_0(\bar x,\bar w)v'_1,\nabla^2_{wx} f_0(\bar x,\bar w)v'_1)\right\}.$$
	Then, $\widehat\Gamma(0)=\{0\}$ if and only if \eqref{Condi_tildeC2} holds. Thus, condition \eqref{Condi_tildeC2} is necessary for $S$ being locally Lipschitz-like around $(\bar w,\bar x)$.

	Let us consider an illustrative example, where the objective function is bilinear and the inequality constraint is polynomial.
	
	\begin{Example} {\rm 
		Consider the problem $(P_w)$ with
		$f_0(x,w)=\displaystyle\sum_{i=1}^{n}x_iw_i$ and $$F(x,w)=\sum_{i=1}^{n}(1-w_i^2)x_i^2-1$$ for all $(x,w)\in \R^n\times\R^n.$
		The stationary point set of this problem is given by
		\begin{equation}\label{KKT_eq_Ex1}
		S(w)=\{x\in\R^n\mid 0\in \nabla_{x}f_0(x,w)+N_{C(w)}(x)\}=\{x\in\R^n\mid -w\in N_{C(w)}(x)\}.
		\end{equation}
		In particular, for $\bar w=0$ one has $S(\bar w)=\{x\in\R^n\mid \|x\|\leq 1\}$. Let $\bar x\in S(\bar x)$ and $\|\bar x\|<1$. Note that
		$$\nabla^{2}_{xx}f_0(\bar x,\bar w)=0,\ \,\nabla^{2}_{wx}f_0(\bar x,\bar w)=E,$$ where $E$ stands for the unit matrix in $\R^{n\times n}$. Since
		${\rm ker}\,\nabla^{2}_{xx}f_0(\bar x,\bar w)=\R^n$, the sufficient condition \eqref{Aubin_Suffi_Condi1} fails. However, we can assert that $S$ is not locally Lipschitz-like around $(\bar w,\bar x)$ because  the necessary condition \eqref{Condi_tildeC2} is not satisfied. In fact, we  directly prove that $S$ is not locally Lipschitz-like around $(\bar w,\bar x)$. Indeed, from \eqref{KKT_eq_Ex1} there exist  neighborhoods $W$ of $\bar w$ and $U$ of $\bar x$ with $S(w)\cap V=\emptyset$ for all $w\in W\setminus\{\bar w\}$. This leads us to the desired result.}
	\end{Example}
	
	\begin{Remark} {\rm The above arguments show that if $(\bar x,\bar w)\in{\rm int}\,D$, then
		\begin{equation}\label{Regularity_M}
		D^*M(\bar\tau)(v'_1)=\widehat D^*M(\bar\tau)(v'_1)=\{0\}
		\end{equation} for any $v'_1\in \R^n$. This implies that $M$ is graphically regular at $\bar\tau$. According to Lemma~ \ref{Combined_estimates}, if the constraint qualification \eqref{C1_condition} is valid, then $\widehat\Gamma(x')=D^*S(\bar w|\bar x)(x')=\Gamma(x')$
		for any $x'\in \R^n$. Suppose that \eqref{C1_condition} is fulfilled. Then,
		$$D^*S(\bar w|\bar x)(x')=\displaystyle\bigcup_{v'_1\in \R^n}\left\{w'\in \R^d\mid (-x',w')\in \nabla G(\bar x,\bar w)^*v'_1\right\}$$ for any $x'\in \R^n$.
		In particular,
		$$D^*S(\bar w|\bar x)(0)=
		\left\{w'\in \R^d\mid\exists v'_1\in \R^n\ {\rm with}\ (0,w')=\nabla G(\bar x,\bar w)^*v'_1\right\}.$$
		According to the Mordukhovich criterion, $S$ is locally Lipschitz-like around $(\bar w,\bar x)$ if and only if $D^*S(\bar w|\bar x)(0)=\{0\}$. The latter means that
		$$\left\{w'\in \R^d\mid \exists v'_1\in \R^n\ {\rm with}\ \nabla^{2}_{xx}f_0(\bar x,\bar w)v'_1=0,\; w'=\nabla^{2}_{wx}f_0(\bar x,\bar w)v'_1\right\}=\{0\}.$$
		Clearly, this is fulfilled if and only if  \eqref{Condi_tildeC2} is satisfied. Now let us verify the constraint qualification \eqref{C1_condition}. Due to  \eqref{Regularity_M},  \eqref{C1_condition} becomes
		$$0\in \nabla G(\bar x,\bar w)^*v'_1\ \;\Longrightarrow\ \;v'_1=0,$$
		or, equivalently,
		$$0\in \nabla^{2}f_0(\bar x,\bar w)^*(v'_1,0)\ \; \Longrightarrow\ \;v'_1=0.$$
		The last condition has been proved to be equivalent to \eqref{Condi_tildeC1_S}. Thus, if \eqref{Condi_tildeC1_S} is satisfied, then  $S$ is locally Lipschitz-like around $(\bar w,\bar x)$ if and only if \eqref{Condi_tildeC2} holds.}
	\end{Remark}
	
	The following theorem summarizes our results for the case of interior points.
	
	\begin{Theorem}\label{Criteria_No1} Suppose that $F(\bar x,\bar w)<0$. The following assertions are valid:
		\begin{description}
			\item[{\rm (a)}] If $S$ is locally Lipschitz-like around $(\bar w,\bar x)$, then condition \eqref{Condi_tildeC2} holds;
			\item[{\rm (b)}] If conditions \eqref{Condi_tildeC1_S} and \eqref{Condi_tildeC2} are simultaneously fulfilled, then  $S$ is locally Lipschitz-like around $(\bar w,\bar x)$;
			\item[{\rm (c)}] If condition \eqref{Condi_tildeC1_S} is satisfied, then  $S$ is locally Lipschitz-like around $(\bar w,\bar x)$ if and only if condition \eqref{Condi_tildeC2} holds.
		\end{description}
	\end{Theorem}
	
	Thus, if condition \eqref{Condi_tildeC1_S} fails we can assert nothing about the local Lipschitz-likeness of $S$ around $(\bar w,\bar x)$.  The next example shows that $S$ can be locally Lipschitz-like around $(\bar w,\bar x)$  when \eqref{Condi_tildeC1_S} is not satisfied,
	
	\begin{Example}  {\rm 
		Consider $(P_w)$ with
		$f_0(x,w)=\frac{1}{3}x^3-w^2x$ and $$F(x,w)=x^2+w^2-1$$ for all $(x,w)\in \R\times\R$. Put $\bar w=0$. The point $\bar x=0$ belongs to $S(\bar w)$ because $(\bar x,\bar w)\in{\rm int}\,D$ and $\nabla_xf_0(\bar x,\bar w)=0$. Since $\nabla^{2}_{xx}f_0(\bar x,\bar w)=0$ and $$\nabla^{2}_{wx}f_0(\bar x,\bar w)=0,$$ condition \eqref{Condi_tildeC1_S} is invalid. We have known that $S$ is locally defined by
		$$\begin{array}{rl}
		S(w)&=\{x\in\R\mid \nabla_xf_0(x,w)=0\}.\\
		&=\{x\in\R\mid |x|=|w|\}.
		\end{array}$$
		Then, for any $x'\in\R$,
		$$\begin{array}{rl}
		D^*S(\bar w|\bar x)(x')&=\{w'\in\R\mid (w',-x')\in N_{{\rm gph}\,S}(0,0)\}.\\
		& =\{w'\in\R\mid |w'|=|x'|\}.
		\end{array}$$
		It follows that $D^*S(\bar w|\bar x)(0)=\{0\}$. Thus, $S$ is locally Lipschitz-like around $(\bar w,\bar x)$.}
	\end{Example}
		
		One referee of the present paper asks: \textit{Whether the condition \eqref{Condi_tildeC2} alone can ensure the local Lipschitz-likeness of $S$ around $(\bar w,\bar x)$.} Answering the question, we construct the next example to demonstrate that, even for polynomial optimization problems, \eqref{Condi_tildeC2} is not sufficient for the later property to hold.
	
	\begin{Example}  {\rm 
		Consider $(P_w)$ with
		$f_0(x,w)=\frac{1}{4}x^4-w^2x$ and $$F(x,w)=x-w-1$$ for all $(x,w)\in \R\times\R$. Put $\bar w=0$. The point $\bar x:=0$ belongs to $S(\bar w)$ because $F(\bar x,\bar w)<0$ and $\nabla_xf_0(\bar x,\bar w)=0$. Since $\nabla^{2}_{xx}f_0(\bar x,\bar w)=\nabla^{2}_{wx}f_0(\bar x,\bar w)=0$, we see that \eqref{Condi_tildeC1_S} fails, but \eqref{Condi_tildeC2} is valid. As $S(w)=\{x\in\R\mid x^3-w^2=0\}=\{w^{\frac{2}{3}}\}$, the stationary point set map is not locally Lipschitz-like around $(\bar w,\bar x)=(0,0)$.}
	\end{Example}
	
	\begin{Remark} {\rm The second assertion of Theorem \ref{Criteria_No1} can be obtained by the classical implicit function theorem. Indeed, if \eqref{Condi_tildeC1_S} and \eqref{Condi_tildeC2} are simultaneously fulfilled, i.e., \eqref{Aubin_Suffi_Condi1} is satisfied, then by \cite[Theorem 1B.1]{DR_2009} (see also \cite[Theorem~9.28]{Rudin_1976}) the implicit multifunction $w\mapsto \{x\in\R^n\mid \nabla_xf_0(x,w)=0\}$ defined by the equation $\nabla_xf_0(x,w)=0$ has a \textit{single-valued localization}  \cite[p.~4]{DR_2009} around $\bar w$ for $\bar x$ which is continuously differentiable in a neighborhood of $\bar w$. This means that there exist a neighborhood $W$ of $\bar w$ and a neighborhood $U$ of $\bar x$ such that for each $w\in W$ there is a unique vector $x=s(w)$ in $U$ satisfying the equation $\nabla_xf_0(x,w)=0$ and $s:W\to U$ is continuously differentiable. Without loss of generality, we can assume that $F(x,w) <0$ for all $(x,w)\in U\times W$. So, by \eqref{new_S}, $S(w)\cap U=\{s(w)\}$ for all $w\in W$. Hence, $S$ is locally Lipschitz-like around $(\bar w,\bar x)$.}
	\end{Remark}
	
	\begin{Remark} {\rm Consider the extended stationary point set map $(w,z)\mapsto \widetilde S(w,z)$ of $(P_w)$, which is defined by $$\widetilde S(w,z)= \{x\in \R^n\mid z\in \nabla_{x}f_0(x,w)+N_{C(w)}(x)\}.$$ If $F(\bar x,\bar w)<0$, then around the point  $((\bar w,0),\bar x )\in {\rm gph}\,\widetilde S$ one can represent $\bar S$ locally as 
			 $$\widetilde S(w,z)= \{x\in \R^n\mid \widetilde G(x,w,z)=0\},$$ where $\widetilde G(x,w,z):=\nabla_{x}f_0(x,w)-z$. Since $\nabla_{(w,z)}\widetilde G(\bar x,\bar w,0)$ has full rank, by \cite[Theorem~2.1]{LeMo04} one has 
			 	$$\begin{array}{rcl}D^*\widetilde S((\bar w,0),\bar x )(x')=\{(w',z')\in \R^d\times\R^n & \mid & w'=-\big(\nabla_{wx}^2f_0(\bar x,\bar w)\big)^Tz',\\
			 	& & x'=\nabla_{xx}^2f_0(\bar x,\bar w)z'\}.\end{array}$$ Hence,
			 	$$\begin{array}{rcl}D^*\widetilde S((\bar w,0),\bar x)(0)=\{(w',z')\in \R^d\times\R^n & \mid & \nabla_{xx}^2f_0(\bar x,\bar w)z'=0\\
			 	& & w'=-\big(\nabla_{wx}^2f_0(\bar x,\bar w)\big)^Tz'\}.\end{array}$$
			Therefore, $D^*\widetilde S((\bar w,0),\bar x)(0)=\{0\}$ if and only if $
			{\rm ker}\,\nabla^{2}_{xx}f_0(\bar x,\bar w)=\{0\}.$ Thanks to Mordukhovich's criterion, this shows that 	
			condition \eqref{Aubin_Suffi_Condi1} is necessary and sufficient for the local Lipschitz-like
			property of the the extended stationary point set map $\widetilde S$	
			around $((\bar w,0),\bar x )$.}
	\end{Remark}
	
\section{Boundary points}

Suppose that $F(\bar x,\bar w)=0$, i.e., $(\bar x,\bar w)$ is a boundary point of $\mathcal{D}$. To obtain a sufficient condition for the Lipschitz-like property of $S$ around $(\bar w,\bar x)$, we will follow the scheme which has been used in the case of interior points. We first need to find out a condition which guarantees the fulfillment of \eqref{C0_condition}. Note that \textbf{(MFCQ)} yields  $\nabla_{x}F(\bar x,\bar w)\neq 0$. Since $\partial g(F(\bar x,\bar w))\subset\R$ and $$\nabla_{x}F(\bar x,\bar w)^*\gamma=\gamma\nabla_{x}F(\bar x,\bar w)$$ for any $\gamma\in\R$, there exists only one element $\lambda\in \partial g(F(\bar x,\bar w))$ satisfying
\begin{equation}\label{Lagrange_multiplier}
\nabla_{x}f_0(\bar x,\bar w)+\lambda\nabla_{x}F(\bar x,\bar w)=0.
\end{equation} This element $\lambda$ is the unique \textit{Lagrange multiplier} for the stationary point $\bar x$ of the minimization problem $\big(P_{\bar w}\big)$. Due to $\lambda\in \partial g(F(\bar x,\bar w))$, we have $\lambda\geq 0$. Let us consider two possibilities of the Lagrange multiplier.

\medskip
\noindent {\bf 4.1 The nondegenerate case}

\smallskip
\noindent Suppose that $\lambda>0$.  Clearly, the equality ${\rm gph}\,\partial g=(\R_{-}\times \{0\})\cup (\{0\}\times\R_{+})$ yields
$$N_{{\rm gph}\,\partial g}(F(\bar x,\bar w),\lambda)=N_{{\rm gph}\,\partial g}(0,\lambda)=\R\times \{0\}.$$
By definition, for any $v'\in\R^n\times\R^d$,
\begin{eqnarray*}\begin{array}{rl} & \partial^2g(F(\bar x,\bar w)|\,\lambda)(\nabla F(\bar x,\bar w)v')\\ & =D^*(\partial g)(F(\bar x,\bar w)|\lambda)(\nabla F(\bar x,\bar w)v')\\
		& =\{u'\in \R\mid (u',-\nabla F(\bar x,\bar w)v')\in N_{{\rm gph}\,\partial g}(F(\bar x,\bar w),\lambda)\}.
	\end{array}\end{eqnarray*}
	Hence
	\begin{equation*}
	\partial^2g(F(\bar x,\bar w)|\,\lambda)(\nabla F(\bar x,\bar w)v')=\begin{cases}
	\R& \mbox{ if } \nabla F(\bar x,\bar w)v'=0\\\
	\emptyset &  \mbox{ if } \nabla F(\bar x,\bar w)v'\neq 0.
	\end{cases}
	\end{equation*}
	So,
	\begin{equation*}
	\Omega_1(\lambda,v')=\begin{cases}
	\{\gamma\nabla F(\bar x,\bar w)\mid \gamma\in\R\}& \mbox{ if } \nabla F(\bar x,\bar w)v'=0\\
	\emptyset &  \mbox{ if } \nabla F(\bar x,\bar w)v'\neq 0.
	\end{cases}
	\end{equation*}
	For $\bar v:=\nabla F(\bar x,\bar w)^*\lambda$, the conditions $\bar y\in \partial g(F(\bar x,\bar w))$ and  $\nabla F(\bar x,\bar w)^*\bar y=\bar v$ force $\bar y=\lambda$. So, by \eqref{Upper_esti_2nd_order_subdifferential_of_f} we get
	\begin{eqnarray*}\begin{array}{rl} \partial^{2}f((\bar x,\bar w)|\bar v)(v')&\subset \lambda\nabla^2F(\bar x,\bar w)v'+\Omega_1(\lambda,v')\\
			& =	\lambda\nabla^2F(\bar x,\bar w)v'+\{\gamma\nabla F(\bar x,\bar w)\mid \gamma\in\R\}
		\end{array}\end{eqnarray*}
		if $\nabla F(\bar x,\bar w)v'=0$, and $\partial^{2}f((\bar x,\bar w)|\bar v)(v')=\emptyset$ if $\nabla F(\bar x,\bar w)v'\neq 0.$  Since $\bar y=\lambda$ is the unique element satisfying the conditions  $\bar y\in \partial g(F(\bar x,\bar w))$ and  $$\nabla F(\bar x,\bar w)^*\bar y=-\nabla_{x}f_0(\bar x,\bar w),$$ from~\eqref{Omega_2} it follows that
		\begin{equation}\label{Upper_est_for_Omega2}\Omega_2(v'_1)\subset \left\{\lambda\nabla^2F(\bar x,\bar w)(v'_1,0)+\gamma\nabla F(\bar x,\bar w)\mid \gamma\in\R\right\},\end{equation}
		for any $v'_1$ with $\nabla_{x}F(\bar x,\bar w)v'_1=0$. Besides, if $\nabla_{x}F(\bar x,\bar w)v'_1\neq 0$, then the set $\Omega_2(v'_1)$ is empty. Thus, \eqref{C0_condition} is fulfilled if the following is satisfied: \textit{for any} $v'_1\in\R^n$,
		\textit{if}
		\begin{equation*}
		\begin{cases}-\nabla^{2}f_0(\bar x,\bar w)^*(v'_1,0)
		\in\left\{\lambda\nabla^2F(\bar x,\bar w)(v'_1,0)+\gamma\nabla F(\bar x,\bar w)\mid \gamma\in\R\right\}\\
		\nabla_xF(\bar x,\bar w)v'_1=0,
		\end{cases}
		\end{equation*}\textit{then} $v'_1=0$.
		Equivalently, \textit{for any} $v'_1\in\R^n$,
		\textit{if}
		\begin{equation*}%%%\label{}
		\begin{cases}0=\nabla^{2}_{xx}f_0(\bar x,\bar w)v'_1+\lambda\nabla^2_{xx}F(\bar x,\bar w)v'_1+\gamma\nabla_{x}F(\bar x,\bar w),\\
		0=\nabla^{2}_{wx}f_0(\bar x,\bar w)v'_1+\lambda\nabla^2_{wx}F(\bar x,\bar w)v'_1+\gamma\nabla_{w}F(\bar x,\bar w),\\
		\nabla_xF(\bar x,\bar w)v'_1=0,\ \gamma\in\R,
		\end{cases}
		\end{equation*}\textit{then} $v'_1=0$.
		Putting
		$$A_1=\left[\begin{array}{ccc} \nabla^{2}_{xx}f_0(\bar x,\bar w)+\lambda\nabla^{2}_{xx}F(\bar x,\bar w) & \nabla_{x}F(\bar x,\bar w) \end{array}\right]\in \R^{n\times (n+1)}$$
		and
		$$A_2=\left[\begin{array}{ccc} \nabla^{2}_{wx}f_0(\bar x,\bar w)+\lambda\nabla^{2}_{wx}F(\bar x,\bar w) & \nabla_{w}F(\bar x,\bar w)\end{array}\right]\in \R^{d\times (n+1)},$$
		where $\nabla_{x}F(\bar x,\bar w)$ and  $\nabla_{w}F(\bar x,\bar w)$ are interpreted as column vectors, we can rewrite the last condition equivalently as follows:
		\begin{equation*}
		\begin{cases}
		A_1\left(\begin{matrix}
		v'_1\\\gamma
		\end{matrix}\right)=0,\quad
		A_2\left(\begin{matrix}
		v'_1\\\gamma
		\end{matrix}\right)=0\\
		\nabla_xF(\bar x,\bar w)v'_1=0,\quad
		\gamma\in\R
		\end{cases} \quad\quad\Longrightarrow\ \;  v'_1=0.
		\end{equation*}
		Since $\nabla_{x}F(\bar x,\bar w)\neq 0$, the latter is equivalent to saying that
		\begin{equation}\label{Condi_tildeC1_S_haty>0}
		{\rm ker}\,A_1\cap {\rm ker}\,A_2\cap \left({\rm ker}\,\nabla_{x}F(\bar x,\bar w)\times \R\right)=\{(0,0)\}.
		\end{equation}
		We now suppose that condition \eqref{Condi_tildeC1_S_haty>0} is satisfied. Then, by Lemma~\ref{Corollary3.1_LeMo04}, for any $x'\in\R^n$,
		$$D^*S(\bar w|\bar x)(x')\subset \displaystyle\bigcup_{v'_1\in \R^n}\left\{w'\in \R^d\mid (-x',w')-\nabla^{2}f_0(\bar x,\bar w)^*(v'_1,0)\in \Omega_2(v'_1)\right\},$$
		with $\Omega_2(v'_1)$ admitting the upper estimation \eqref{Upper_est_for_Omega2}. So, for any $x'\in\R^n$,
		\begin{equation}\label{Upper_esti_for_S_nc} D^*S(\bar w|\bar x)(x')\subset\Gamma_2(x'),\end{equation}
		where $\Gamma_2(x')$ consists of vectors $w'\in\R^d$ such that
		\begin{equation*}
		\begin{cases}
		w'=\left[\nabla^{2}_{wx}f_0(\bar x,\bar w)+\lambda\nabla^2_{wx}F(\bar x,\bar w)\right]v'_1+\gamma\nabla_{w}F(\bar x,\bar w),\\
		\left[\nabla^{2}_{xx}f_0(\bar x,\bar w)+\lambda\nabla^2_{xx}F(\bar x,\bar w)\right]v'_1+\gamma\nabla_{x}F(\bar x,\bar w)=-x',\\
		\nabla_xF(\bar x,\bar w)v'_1=0,\ v'_1\in \R^n,\ \gamma\in\R.
		\end{cases}
		\end{equation*}
		
		To obtain a lower estimate for the Fr\' echet coderivative values of $S$, we will use some results of Qui \cite{Qui_JOTA2014}. For any $v'_1\in\R^n$ satisfying $\nabla_x F(\bar x,\bar w)v'_1=0$, the arguments given in \cite[pp.~410--412]{Qui_JOTA2014} provide us with the inclusion
		\begin{equation}\label{Inclusion_HatOmegaM} \widehat\Omega_M(v'_1)\subset\widehat D^*M(\bar\tau)(v'_1),	\end{equation}
		where $M$ is the multifunction in \eqref{Sum_for_S} and
		\begin{equation*}
		\begin{array}{rl}
		\widehat\Omega_M(v'_1):= & \{(x',w')\in\R^n\times\R^d\mid x'=\gamma\nabla_{x}F(\bar x,\bar w)+\hat{y} \nabla^{2}_{xx}F(\bar x,\bar w)v'_1,\\
		& \qquad \qquad \quad\ \; w'=\gamma\nabla_{w}F(\bar x,\bar w)+\hat{y}\nabla^{2}_{wx}F(\bar x,\bar w)v'_1,\, \gamma\in\R\}.
		\end{array}
		\end{equation*}
		In addition, if $\nabla_x F(\bar x,\bar w)v'_1\neq 0$, then by the reasoning in \cite[pp.~405--406]{Qui_JOTA2014} we have $\widehat D^*M(\bar\tau)(v'_1)=\emptyset$. Note that the Lagrange multiplier $\lambda$ here coincides with the constant $\mu$ in \cite[Theorem~3.2]{Qui_JOTA2014} and \textit{the above assertions do not require $\nabla_wF(\bar x,\bar w)\neq 0$.} Then, by~\eqref{Widehat_Gamma(x')} and \eqref{Inclusion_HatOmegaM},
		$$\widehat\Gamma(x')\supset\displaystyle\bigcup_{\begin{subarray}{c}v'_1\in \R^n,\\\nabla_xF(\bar x,\bar w)v'_1=0\end{subarray}}\left\{w'\in \R^d\mid (-x',w')\in \nabla G(\bar x,\bar w)^*v'_1+\widehat\Omega_M(v'_1)\right\}.$$
		Since $\nabla G(\bar x,\bar w)^*v'_1=(\nabla^2_{xx} f_0(\bar x,\bar w)v'_1,\nabla^2_{wx} f_0(\bar x,\bar w)v'_1)$, one can easily show that the right-hand-side set equals to $\Gamma_2(x')$. Therefore, if \eqref{Condi_tildeC1_S_haty>0} is satisfied, then by \eqref{Upper_esti_for_S_nc} we have
		$$\widehat\Gamma(x')=\widehat D^*S(\bar w|\bar x)(x')=D^*S(\bar w|\bar x)(x')=\Gamma_2(x').$$
		Thus, under the assumption \eqref{Condi_tildeC1_S_haty>0}, $S$ is locally Lipschitz-like around $(\bar w,\bar x)$ if and only if $\Gamma_2(0)=\{0\}$. Clearly, the set $\Gamma_2(0)$ consists of vectors $w'\in\R^d$ such that
		\begin{equation*}%%%\label{}
		\begin{cases}
		w'=\left[\nabla^{2}_{xw}f_0(\bar x,\bar w)^{T}+\lambda\nabla^2_{wx}F(\bar x,\bar w)\right]v'_1+\gamma\nabla_{w}F(\bar x,\bar w),\\
		\left[\nabla^{2}_{xx}f_0(\bar x,\bar w)+\lambda\nabla^2_{xx}F(\bar x,\bar w)\right]v'_1+\gamma\nabla_{x}F(\bar x,\bar w)=0,\\
		\nabla_xF(\bar x,\bar w)v'_1=0,\ v'_1\in \R^n,\ \gamma\in\R.
		\end{cases}
		\end{equation*}
		Equivalently, $w'\in\R^d$ belongs to $\Gamma_2(0)$ iff
		\begin{equation*}
		\begin{cases}
		w'=A_2\left(\begin{matrix}
		v'_1\\\gamma
		\end{matrix}\right),\\
		A_1\left(\begin{matrix}
		v'_1\\\gamma
		\end{matrix}\right)=0,\ \,\nabla_xF(\bar x,\bar w)v'_1=0,\\
		v'_1\in \R^n,\ \gamma\in\R.
		\end{cases}
		\end{equation*} So, $\Gamma_2(0)=\{0\}$ if and only if
		$A_2\left(\begin{matrix}
		v'_1\\\gamma
		\end{matrix}\right)=0$
		for any pair $(v'_1,\gamma)\in\R^n\times\R$ satisfying
		\begin{equation*}
		A_1\left(\begin{matrix}
		v'_1\\\gamma
		\end{matrix}\right)=0,\ \;\nabla_xF(\bar x,\bar w)v'_1=0.
		\end{equation*}
		The latter happens if and only if
		\begin{equation}\label{Condi_tildeC2_haty>0}
		{\rm ker}\,A_1\cap \left({\rm ker}\,\nabla_{x}F(\bar x,\bar w)\times \R\right)\subset {\rm ker}\,A_2.
		\end{equation}
		
		To sum up, we state the following theorem.
		
		\begin{Theorem}\label{Criteria_No2} Suppose that $F(\bar x,\bar w)=0$ and the Lagrange multiplier $\lambda$ corresponding to the stationary point $\bar x\in S(\bar w)$ is positive. If \eqref{Condi_tildeC1_S_haty>0} holds, then $S$ is locally Lipschitz-like around $(\bar w,\bar x)$ if and only if \eqref{Condi_tildeC2_haty>0} is satisfied.
		\end{Theorem}
		
		\begin{Remark} {\rm 
			Combining the conditions \eqref{Condi_tildeC1_S_haty>0} and \eqref{Condi_tildeC2_haty>0}, we obtain a sufficient condition for $S$ being locally Lipschitz-like around $(\bar w,\bar x)$, that is
			\begin{equation}\label{Aubin_Suffi_Condi2}
			{\rm ker}\,A_1\cap \left({\rm ker}\,\nabla_{x}F(\bar x,\bar w)\times \R\right)=\{0\}.
			\end{equation}}
		\end{Remark}
		
		\begin{Example} {\rm 
			Consider the problem $(P_w)$ with $f_0(x,w)=-x^2+(w-1)x$ and $F(x,w)=x^2+w^2-2$ for any $(x,w)\in\R\times\R$. Then, by \eqref{KKT_point-set}, the stationary point set map of $(P_w)$ is defined by
			$$S(w)=\left\{x\in\R\mid -2x+w-1+\partial_xf(x,w)\right\},$$
			with $f(x,w)=(g\circ F)(x,w)$ and $g(y)=\delta_{\R_{-}}(y)$ for any $y\in\R$. Let $\bar w=1$ and $\bar x=1$. Since $F(\bar x,\bar w)=0$ and $\nabla_xF(\bar x,\bar w)=2$, condition \eqref{MFCQ_condition} is valid. Hence, from \eqref{Implication_of_s_amenability2} we have
			\begin{equation*}
			\begin{array}{rl}
			\partial_xf(\bar x,\bar w)& = \nabla_xF(\bar x,\bar w)^*N_\mathbb{R_{-}}(F(\bar x,\bar w))\\
			& =\nabla_xF(\bar x,\bar w)^*\mathbb{R_{+}}=\mathbb{R_{+}}.
			\end{array}
			\end{equation*}
			Now, it easy to show that $\bar x\in S(\bar w)$. We have $\nabla_xf_0(\bar x,\bar w)=-2$ and $\lambda=1$ due to \eqref{Lagrange_multiplier}. Hence, $A_1=[\begin{matrix} 0 & 2\end{matrix}]$ and ${\rm ker}\,A_1=\R\times \{0\}$. Thus, \eqref{Aubin_Suffi_Condi2} is fulfilled and consequently $S$ is locally Lipschitz-like around $(\bar w,\bar x)$.}
		\end{Example}
		
		\begin{Remark}\label{remark_about_S(w)} {\rm For any stationary point $\bar x\in S(\bar w)$ satisfying \textbf{(MFCQ)}, the corresponding unique multiplier $\lambda$ is defined by the equation \eqref{Lagrange_multiplier}. This fact justifies to the following assertion: \textit{Given any $w\in\R^d$, if $\nabla_xF(x,w)\neq 0$ for all $x$ with $F(x,w)=0$, then one has}
			\begin{equation*}\label{S(w)}
			\begin{array}{rl}
			S(w)  = \big\{x\in {\R}^n\mid  & F(x,w)\leq 0,\ \exists\lambda\geq 0\ {\rm s.t.}\  \lambda F(x,w)=0,\\
			&  \nabla_{x}f_0(x,w)+\lambda\nabla_{x}F(x,w)=0\big\}.
			\end{array}
			\end{equation*}}		
		\end{Remark}
		
		One referee of the present paper asks: \textit{Whether the condition \eqref{Condi_tildeC2_haty>0} alone can ensure the local Lipschitz-likeness of $S$ around $(\bar w,\bar x)$.} To answer the question, let us consider the next example showing that, even for polynomial optimization problems, \eqref{Condi_tildeC2_haty>0}  is not sufficient for $S$ to be locally Lipschitz-like around $(\bar w,\bar x)$.
		
		\begin{Example} {\rm 
			Consider $(P_w)$ with $n=2$, $d=1$, $f_0(x,w)=\frac{1}{4}wx_1^4-wx_1-x_2$, and $F(x,w)=wx_1+x_2-w$ for all $(x,w)=(x_1,x_2,w)\in\R^2\times\R$. Choose $\bar x=(0,0)$ and $\bar w=0$. Using Remark 4.4, one can show that $S(\bar w)= \R\times \{0\}$ and $S(w)= \{(0,w)\}$ for every $w\neq 0$. Hence the stationary map $S(.)$ is not locally Lipschitz-like around $(\bar w,\bar x)$. Furthermore, since the unique Lagrange multiplier corresponding to $\bar x\in S(\bar w)$ is $\lambda=1$, one can find that  ${\rm ker}\, A_1=\R\times \R\times\{0\}$, ${\rm ker}\, A_2=\R\times \R\times \{0\}$, and  ${\rm ker}\,\nabla_xF(\bar w,\bar x)=\R\times\{0\}$. So \eqref{Condi_tildeC2_haty>0} is satisfied, but \eqref{Condi_tildeC1_S_haty>0} fails to hold.}
		\end{Example}
				
		\noindent {\bf 4.2 The degenerate case}
		
		\smallskip
		\noindent
		Consider the second possibility: $\lambda=0$. We have
		$$N_{{\rm gph}\,\partial g}(F(\bar x,\bar w),\lambda)=N_{{\rm gph}\,\partial g}(0,0)=(\{0\}\times \R)\cup (\R\times \{0\})\cup (\R_{+}\times \R_{-}).$$
		Hence, for any $v'\in\R^n\times\R^d$,
		\begin{equation*}
		\partial^2g(F(\bar x,\bar w)|\,\lambda)(\nabla F(\bar x,\bar w)v')=\begin{cases}
		\R& \mbox{ if } \nabla F(\bar x,\bar w)v'=0,\\
		\R_{+}& \mbox{ if } \nabla F(\bar x,\bar w)v'>0,\\
		\{0\} &  \mbox{ if } \nabla F(\bar x,\bar w)v'<0.
		\end{cases}
		\end{equation*}
		Consequently,
		\begin{equation*}
		\Omega_1(\lambda,v')=
		\begin{cases}
		\{\gamma\nabla F(\bar x,\bar w)\mid \gamma\in\R\}& \mbox{ if } \nabla F(\bar x,\bar w)v'=0,\\
		\{\gamma\nabla F(\bar x,\bar w)\mid \gamma\in\mathbb{R_{+}}\}& \mbox{ if } \nabla F(\bar x,\bar w)v'>0,\\
		\{0\} &  \mbox{ if } \nabla F(\bar x,\bar w)v'<0.
		\end{cases}
		\end{equation*}
		In this case, $\bar y=\lambda=0$ is the unique element satisfying $\bar y\in \partial g(F(\bar x,\bar w))$ and  $\nabla_xF(\bar x,\bar w)^*\bar y=-\nabla_{x}f_0(\bar x,\bar w)$. So, by \eqref{Omega_2} we have
		$$\Omega_2(v'_1)=\partial^2f((\bar x,\bar w)|0)(v'_1,0).$$
		Clearly, the conditions $\bar y\in \partial g(F(\bar x,\bar w))$ and $\nabla F(\bar x,\bar w)^*\bar y=0$ imply $\bar y=\lambda=0$. So, for $\bar v=0$, from \eqref{Upper_esti_2nd_order_subdifferential_of_f} we get
		$$\partial^2f((\bar x,\bar w)|0)(v')\subset \Omega_1(\lambda,v')$$ for all $v'\in \R^n\times \R^d$. It follows that
		\begin{equation}\label{Omega2_degenerate}
		\Omega_2(v'_1)\subset
		\begin{cases}
		\{\gamma\nabla F(\bar x,\bar w)\mid \gamma\in\R\}& \mbox{ if } \nabla_{x}F(\bar x,\bar w)v'_1=0,\\
		\{\gamma\nabla F(\bar x,\bar w)\mid \gamma\in\mathbb{R_{+}}\}& \mbox{ if } \nabla_{x}F(\bar x,\bar w)v'_1>0,\\
		\{0\} &  \mbox{ if } \nabla_{x}F(\bar x,\bar w)v'_1<0.
		\end{cases}
		\end{equation}
		Hence, the condition \eqref{C0_condition} is satisfied if the following holds:
		\textit{For any} $v'_1\in \R^n$,
		\begin{description}
			\item[{\rm (i)}] \textit{if} $\nabla_{x}F(\bar x,\bar w)v'_1=0$ \textit{and}
			\begin{equation*}
			\left(-\nabla^{2}_{xx}f_0(\bar x,\bar w)v'_1,-\nabla^{2}_{wx}f_0(\bar x,\bar w)v'_1\right)\in \{\gamma\nabla F(\bar x,\bar w)\mid \gamma\in\R\},
			\end{equation*} \textit{then} $v'_1=0$;
			\item[{\rm (ii)}] \textit{if} $\nabla_{x}F(\bar x,\bar w)v'_1>0$, \textit{then} $$\left(-\nabla^{2}_{xx}f_0(\bar x,\bar w)v'_1,-\nabla^{2}_{wx}f_0(\bar x,\bar w)v'_1\right)\notin \{\gamma\nabla F(\bar x,\bar w)\mid \gamma\in\mathbb{R_+}\};$$
			\item[{\rm (iii)}] \textit{if} $ \nabla_{x}F(\bar x,\bar w)v'_1<0$, \textit{then} $\left(\nabla^{2}_{xx}f_0(\bar x,\bar w)v'_1,\nabla^{2}_{wx}f_0(\bar x,\bar w)v'_1\right)\neq 0$.
		\end{description}
		It follows that \eqref{C0_condition} is fulfilled if the following holds: \textit{for any} $v'_1\in\R^n$, \textit{if} $$\left(-\nabla^{2}_{xx}f_0(\bar x,\bar w)v'_1,-\nabla^{2}_{wx}f_0(\bar x,\bar w)v'_1\right)\in \{\gamma\nabla F(\bar x,\bar w)\mid \gamma\in\R\},$$
		\textit{then} $v'_1=0$. The latter can be written equivalently as
		\begin{equation*}
		\begin{cases}
		A'_1\left(\begin{matrix}
		v'_1\\\gamma
		\end{matrix}\right)=0,\quad
		A'_2\left(\begin{matrix}
		v'_1\\\gamma
		\end{matrix}\right)=0\\
		v'_1\in\R^n,\quad \gamma\in\R
		\end{cases} \quad\quad\Longrightarrow\ \;  v'_1=0,
		\end{equation*}
		where
		$$A'_1:=\left[\begin{array}{ccc} \nabla^{2}_{xx}f_0(\bar x,\bar w) & \nabla_{x}F(\bar x,\bar w)  \end{array}\right]\in \R^{n\times (n+1)},$$
		and
		$$A'_2:=\left[\begin{array}{ccc} \nabla^{2}_{wx}f_0(\bar x,\bar w) & \nabla_{w}F(\bar x,\bar w)  \end{array}\right]\in \R^{d\times (n+1)}.$$ Since $\nabla_{x}F(\bar x,\bar w)\neq 0$, this condition is equivalent to
		\begin{equation}\label{Condi_tildeC1_S_haty=0}
		{\rm ker}\,A'_1\cap {\rm ker}\,A'_2=\{0\}.
		\end{equation}
		We now suppose that \eqref{Condi_tildeC1_S_haty=0} is valid. Then, by Lemma~\ref{Corollary3.1_LeMo04} and the upper estimate \eqref{Omega2_degenerate}, for any $x'\in\R^n$,
		\begin{equation}\label{Coder_esti_degenerate} D^*S(\bar w|\bar x)(x')\subset\Gamma_3(x'),\end{equation}
		where $\Gamma_3(x')$ consists of vectors $w'\in\R^d$ for which there exists $v'_1\in\R^n$ with
		\begin{equation*}
		\begin{cases}\left(-x'-\nabla^{2}_{xx}f_0(\bar x,\bar w)v'_1,w'-\nabla^{2}_{wx}f_0(\bar x,\bar w)v'_1\right)\in \{\gamma\nabla F(\bar x,\bar w)\mid \gamma\in\R\},\\
		\nabla_{x}F(\bar x,\bar w)v'_1=0,
		\end{cases}
		\end{equation*} or
		\begin{equation*}
		\begin{cases}\left(-x'-\nabla^{2}_{xx}f_0(\bar x,\bar w)v'_1,w'-\nabla^{2}_{wx}f_0(\bar x,\bar w)v'_1\right)\in \{\gamma\nabla F(\bar x,\bar w)\mid \gamma\in\R_+\},\\
		\nabla_{x}F(\bar x,\bar w)v'_1>0,
		\end{cases}
		\end{equation*} or
		\begin{equation*}
		\begin{cases}
		-x'-\nabla^{2}_{xx}f_0(\bar x,\bar w)v'_1=0,\ \, w'-\nabla^{2}_{wx}f_0(\bar x,\bar w)v'_1=0,\\
		\nabla_{x}F(\bar x,\bar w)v'_1<0.
		\end{cases}
		\end{equation*}
		By putting
		$$\Delta_1=\{(v'_1,\gamma)\in\R^n\times\R\mid\nabla_xF(\bar x,\bar w)v'_1>0,\, \gamma\geq 0\}$$
		and
		$$\Delta_2=\{v'_1\in\R^n\mid\nabla_xF(\bar x,\bar w)v'_1<0\},$$
		we see that $\Gamma_3(x')$ consists of vectors $w'\in\R^d$ such that
		\begin{equation*}%%%\label{}
		\begin{cases}
		w'=\nabla^{2}_{wx}f_0(\bar x,\bar w)v'_1 +\gamma\nabla_{w}F(\bar x,\bar w),\\
		\nabla^{2}_{xx}f_0(\bar x,\bar w)v'_1+\gamma\nabla_{x}F(\bar x,\bar w)=-x',\\
		(v'_1,\gamma)\in {\rm ker}\,\nabla_{x}F(\bar x,\bar w)\times \R,
		\end{cases}
		\end{equation*} or
		\begin{equation*}%%%\label{}
		\begin{cases}
		w'=\nabla^{2}_{wx}f_0(\bar x,\bar w)v'_1 +\gamma\nabla_{w}F(\bar x,\bar w),\\
		\nabla^{2}_{xx}f_0(\bar x,\bar w)v'_1+\gamma\nabla_{x}F(\bar x,\bar w)=-x',\\
		(v'_1,\gamma)\in \Delta_1,
		\end{cases}
		\end{equation*} or
		\begin{equation*}
		\begin{cases}
		w'=\nabla^{2}_{wx}f_0(\bar x,\bar w)v'_1,\\
		\nabla^{2}_{xx}f_0(\bar x,\bar w)v'_1=-x',\\
		v'_1\in\Delta_2.
		\end{cases}
		\end{equation*}
		In particular, the set $\Gamma_3(0)$ consists of vectors $w'\in\R^d$ such that
		\begin{equation*}
		\begin{cases}
		w'=A'_2\left(\begin{matrix}
		v'_1\\\gamma
		\end{matrix}\right),\ \,
		A'_1\left(\begin{matrix}
		v'_1\\\gamma
		\end{matrix}\right)=0,\\
		(v'_1,\gamma)\in {\rm ker}\,\nabla_{x}F(\bar x,\bar w)\times \R.
		\end{cases}
		\end{equation*} or
		\begin{equation*}
		\begin{cases}
		w'=A'_2\left(\begin{matrix}
		v'_1\\\gamma
		\end{matrix}\right),\ \,
		A'_1\left(\begin{matrix}
		v'_1\\\gamma
		\end{matrix}\right)=0,\\
		(v'_1,\gamma)\in\Delta_1.
		\end{cases}
		\end{equation*} or
		\begin{equation*}
		\begin{cases}
		w'=\nabla^{2}_{wx}f_0(\bar x,\bar w)v'_1,\\
		\nabla^{2}_{xx}f_0(\bar x,\bar w)v'_1=0,\,v'_1\in\Delta_2.
		\end{cases}
		\end{equation*}
		Therefore, $\Gamma_3(0)=\{0\}$ if and only if the following three conditions are simultaneously fulfilled
		\begin{equation}\label{CondiC2a_haty=0}
		{\rm ker}\,A'_1\cap \left({\rm ker}\,\nabla_{x}F(\bar x,\bar w)\times \R\right)\subset {\rm ker}\,A'_2,
		\end{equation}
		\begin{equation}\label{CondiC2b_haty=0}
		{\rm ker}\,A'_1\cap\Delta_1\subset {\rm ker}\,A'_2,
		\end{equation} and
		\begin{equation}\label{CondiC2c_haty=0}
		{\rm ker}\,\nabla^{2}_{xx}f_0(\bar x,\bar w)\cap\Delta_2\subset{\rm ker}\,\nabla^{2}_{wx}f_0(\bar x,\bar w).
		\end{equation}
		Remember that under condition \eqref{Condi_tildeC1_S_haty=0} we have the estimate \eqref{Coder_esti_degenerate}. Therefore, the fulfillment of \eqref{Condi_tildeC1_S_haty=0}, \eqref{CondiC2a_haty=0}, \eqref{CondiC2b_haty=0}, and \eqref{CondiC2c_haty=0} implies $D^*S(\bar w|\bar x)(0)=\{0\}$. Applying the Mordukhovich criterion, we obtain the following\textit{ sufficient condition} for the local Lipschitz-likeness of $S$ around $(\bar w,\bar x)$:
		\begin{equation}\label{Aubin_Suffi_Condi3}
		\begin{cases}
		{\rm ker}\,A'_1\cap {\rm ker}\,A'_2=\{0\},\\
		{\rm ker}\,A'_1\cap \left({\rm ker}\,\nabla_{x}F(\bar x,\bar w)\times \R\right)\subset {\rm ker}\,A'_2,\\
		{\rm ker}\,A'_1\cap\Delta_1\subset {\rm ker}\,A'_2,\\
		{\rm ker}\,\nabla^{2}_{xx}f_0(\bar x,\bar w)\cap\Delta_2\subset{\rm ker}\,\nabla^{2}_{wx}f_0(\bar x,\bar w).
		\end{cases}
		\end{equation}

		Now, by the arguments of Qui \cite[pp.~414--416]{Qui_JOTA2014}, the Fr\'echet coderivative values of the multifunction $M$ in \eqref{Sum_for_S} admit the lower estimate
		$$\widehat D^*M(\bar\tau)(v'_1)\supset \left\{\gamma\big(\nabla_xF(\bar x,\bar w),\nabla_wF(\bar x,\bar w)\big)\mid\gamma\geq 0\right\}$$
		for any $v'_1$ with $\nabla_xF(\bar x,\bar w)v'_1\geq 0$. If $\nabla_xF(\bar x,\bar w)v'_1<0$, then $\widehat D^*M(\bar\tau)(v'_1)=\emptyset$; see \cite[p.~412]{Qui_JOTA2014}. (It is important to stress that \textit{we don't need the condition $\nabla_wF(\bar x,\bar w)\neq 0$ here}.) So, by \eqref{Widehat_Gamma(x')} one has
		$\widehat\Gamma(x')\supset \widehat\Gamma_1(x')$ for any $x'\in\R^n$,
		where
		$$\begin{array}{rl}
		\widehat\Gamma_1(x') & :=\displaystyle\bigcup_{\begin{subarray}{c}v'_1\in \R^n,\\ \nabla_xF(\bar x,\bar w)v'_1\geq 0\end{subarray}}\Big\{w'\in \R^d\mid (-x',w')\\
		& \qquad \qquad \qquad \in \nabla G(\bar x,\bar w)^*v'_1+\left\{\gamma\big(\nabla_xF(\bar x,\bar w),\nabla_wF(\bar x,\bar w)\big)\mid\gamma\geq 0\right\}\Big\}.
		\end{array}$$
		Choosing $v'_1=0$, we have $0\in\widehat\Gamma_1(0)$. In combination with the results in Lemma~\ref{Combined_estimates}, this yields
		$$\{0\}\subset\widehat\Gamma_1(0)\subset\widehat\Gamma(0)\subset\widehat D^*S(\bar w|\bar x)(0)\subset D^*S(\bar w|\bar x)(0).$$ According to the Mordukhovich criterion, if $S$ is locally Lipschitz-like around $(\bar w,\bar x)$, then $D^*S(\bar w|\bar x)(0)=\{0\}$ and $\widehat\Gamma_1(0)=\{0\}$ as a result. Since $$\nabla G(\bar x,\bar w)^*v'_1=(\nabla^2_{xx} f_0(\bar x,\bar w)v'_1,\nabla^2_{wx} f_0(\bar x,\bar w)v'_1),$$ the latter means that
		$$\begin{cases}
		w'=\nabla^{2}_{wx}f_0(\bar x,\bar w)v'_1 +\gamma\nabla_{w}F(\bar x,\bar w)\\
		\nabla^{2}_{xx}f_0(\bar x,\bar w)v'_1+\gamma\nabla_{x}F(\bar x,\bar w)=0\qquad\quad\ \Longrightarrow\ \,w'=0.\\
		\nabla_xF(\bar x,\bar w)v'_1\geq 0,\ \, v'_1\in\R^n,\ \, \gamma\in\R_+
		\end{cases}$$
		For $\Delta_3:=\{(v'_1,\gamma)\in\R^n\times\R\mid\nabla_xF(\bar x,\bar w)v'_1\geq 0,\, \gamma\geq 0\}$,
		this condition becomes
		\begin{equation*}
		\begin{cases}
		w'=A'_2\left(\begin{matrix}
		v'_1\\\gamma
		\end{matrix}\right)\\
		A'_1\left(\begin{matrix}
		v'_1\\\gamma
		\end{matrix}\right)=0\qquad\quad\ \Longrightarrow\ \,w'=0.\\
		\left(v'_1,\gamma\right)\in \Delta_3
		\end{cases}
		\end{equation*} Clearly, the last property is equivalent to
		\begin{equation}\label{Necessary_condi_degenerate}
		{\rm ker}\,A'_1\cap\Delta_3\subset {\rm ker}\,A'_2.
		\end{equation}
		Thus, we have shown that \eqref{Necessary_condi_degenerate} is a \textit{necessary condition} for $S$ being locally Lipschitz-like around $(\bar w,\bar x)$.
		
		Condition \eqref{Aubin_Suffi_Condi3} implies \eqref{Necessary_condi_degenerate}. Indeed, suppose that \eqref{Aubin_Suffi_Condi3} is fulfilled and $\left(v'_1,\gamma\right)\in {\rm ker}\,A'_1\cap\Delta_3$. If $\nabla_xF(\bar x,\bar w)v'_1=0$, then 	$$\left(v'_1,\gamma\right)\in {\rm ker}\,A'_1\cap \left({\rm ker}\,\nabla_{x}F(\bar x,\bar w)\times \R\right);$$ so $\left(v'_1,\gamma\right)\in {\rm ker}\,A'_2$ by the second condition in \eqref{Aubin_Suffi_Condi3}. If $\nabla_xF(\bar x,\bar w)v'_1>0$, then $\left(v'_1,\gamma\right)\in 	{\rm ker}\,A'_1\cap\Delta_1$; hence $\left(v'_1,\gamma\right)\in {\rm ker}\,A'_2$ by the third condition in \eqref{Aubin_Suffi_Condi3}. We have thus proved that \eqref{Aubin_Suffi_Condi3} yields \eqref{Necessary_condi_degenerate}.

		The above elementary analysis clearly shows how the sufficient condition for  the local Lipschitz-likeness of $S$ around $(\bar w,\bar x)$ in the degenerate case is stronger than the necessary one.

		We can summarize our results for the degenerate case as follows.
		
		\begin{Theorem}\label{Criteria_No3} Suppose that $F(\bar x,\bar w)=0$ and the Lagrange multiplier $\lambda$ corresponding to the stationary point $\bar x\in S(\bar w)$ is null. The following assertions are true:
			\begin{description}
				\item[{\rm (a)}] If $S$ is locally Lipschitz-like around $(\bar w,\bar x)$, then condition \eqref{Necessary_condi_degenerate} holds;
				\item[{\rm (b)}] If condition \eqref{Aubin_Suffi_Condi3} is fulfilled, then  $S$ is locally Lipschitz-like around $(\bar w,\bar x)$.
			\end{description}
		\end{Theorem}
		
		Let us consider a simple example to see how Theorem~\ref {Criteria_No3} works for concrete optimization problems.
		
		\begin{Example}\label{Exam_nabla_xF(barx,barw)=0} {\rm 
			Let
			$f_0(x,w)=x^2(w-2)$ and $F(x,w)=w(x-1)$ for all $(x,w)\in \R\times\R$.
			The stationary point set of $(P_w)$ is given by
			\begin{equation*}\label{KKT_eq_Ex2}
			S(w)=\{x\in \R\mid 0\in  2x(w-2)+\partial_{x}f(x,w)\},
			\end{equation*}
			with $f(x,w)=(g\circ F)(x,w)$ and  $g(y)=\delta_{\R_-}(y)$ for all $y\in\R$. Let $\bar w=1$. Then, the point $\bar x=1$ belongs to $S(\bar w)$. Indeed, since $F(\bar x,\bar w)=0$ and $\nabla_xF(\bar x,\bar w)=1$, condition \textbf{(MFCQ)} is valid. Hence, from \eqref{Implication_of_s_amenability2} we have
			\begin{equation*}
			\begin{array}{rl }
			\partial_xf(\bar x,\bar w) & = \nabla_xF(\bar x,\bar w)^*N_\mathbb{R_{-}}(F(\bar x,\bar w))\\
			& =\nabla_xF(\bar x,\bar w)^*\mathbb{R_{+}}=\mathbb{R_{+}}.
			\end{array}
			\end{equation*}
			Now it is easy to check that $\bar x\in S(\bar w)$. Here we have
			$A'_1=[\begin{array}{ccc}-4 & 1\end{array}]$,
			$A'_2=[\begin{array}{ccc} 2 & 0\end{array}]$, and
			${\rm ker}\,\nabla_{x}F(\bar x,\bar w)=\{0\}$. Thus, condition \eqref{Aubin_Suffi_Condi3} is fulfilled; as a result, $S$ is locally Lipschitz-like around $(\bar w,\bar x)$.}
		\end{Example}
		
		\begin{Remark} {\rm 
			Looking back to Example \ref{Exam_nabla_xF(barx,barw)=0}, we see that $$\nabla_wF(\bar x,\bar w)=\bar x-1=0.$$ So, the preceding results of Qui \cite[Theorem~4.2]{Qui_JoGO_2016} cannot be applied for the boundary point $(\bar x,\bar w)$ of the set $$\mathcal{D}=\{(x,w)\mid F(x,w)\leq 0\}=\{(x,w)\mid x\leq 1,\, w\geq 0\}\cup \{(x,w)\mid x\leq 1,\, w\geq 0\},$$ as it requires that $\nabla_wF(\bar x,\bar w)\neq 0$.}
		\end{Remark}
 
 \begin{acknowledgements}	
 	This work was supported by National Foundation for Science $\&$ Technology Development (Vietnam) and the Grant MOST 105-2115-M-039-002-MY3 (Taiwan). The authors are grateful to the anonymous referees for their careful readings and valuable suggestions. Examples~3.4 and~4.2 in this paper present our solutions to two open questions raised by one of the referees. In addition, Remarks 3.2 and 3.3 are based on some comments of that referee.
 \end{acknowledgements}

\end{document}